\input amstex
\document
\magnification=\magstep1
\hbox{}

\def\reel{\text{ }\hbox{{\rm R}\kern-.9em\hbox{{\rm I}}\text{ }}\text{ }}
\def\entier{\text{ }\hbox{{\rm N}\kern-.9em\hbox{{\rm I}}\text{ }}\text{ }}
\def\scalaire{\text{ }\hbox{{\rm K}\kern-.9em\hbox{{\rm I}}\text{ }}\text{ }}
\def\rationel{\text{ }\hbox{{\rm Q}\kern-.9em\hbox{{\rm I}}\text{ }}\text{ }}

\def\bn{\bigskip\noindent}

\hbox{}
\vskip1truecm
\centerline{\bf RENORMINGS OF $L_p(L_q)$}
\bigskip
\centerline{by}
\medskip
\centerline{\bf R. Deville$^*$, R. Gonzalo$^{**}$ and J. A. Jaramillo$^{**}$}

\vskip1truecm
\noindent
(*) Laboratoire de Math\'ematiques,
Universit\'e Bordeaux I,
351, cours de la Lib\'eration,
33400 Talence,
FRANCE.

\smallskip\noindent
email : {\it deville\@math.u-bordeaux.fr}

\bigskip\noindent
(**) Departamento de An\'alisis Matem\'atico,
Universidad Complutense de Madrid,
28040 Madrid,
SPAIN.

\smallskip\noindent
email : {\it raquel26\@eucmax.sim.ucm.es}\quad and\quad {\it
jaramil\@eucmax.sim.ucm.es}

\vskip2truecm
\noindent
{\bf Abstract} : {\it We investigate the best order of smoothness of $L_p(L_q)$.
We prove in particular that there exists a $\Cal C^\infty$-smooth bump function
on $L_p(L_q)$ if and only if $p$ and $q$ are even integers and $p$ is a multiple
of $q$. }

\vskip 6.5cm\noindent
1991 \it Mathematics Subject Classification \rm : 46B03, 46B20, 46B25, 46B30.

\noindent
\it Key words \rm : Higher order smoothness, Renormings, Geometry of Banach
spaces, $L_p$ spaces.

\bn
(**) The research of second and the third author has been
partially supported by DGICYT grant PB 93-0452.

\newpage

\noindent
{\bf Introduction }
\medskip
The first results on high order differentiability of
the norm in a Banach space were given by
Kurzweil [10] in 1954, and Bonic and Frampton [2] in 1965. In
[2] the best order of differentiability of an equivalent norm
on classical $L_p$-spaces, $1<p<\infty$, is given.
The best order of smoothness of equivalent norms and of bump functions
in Orlicz spaces has been investigated by Maleev and Troyanski
(see [13]).

The  work of Leonard and Sundaresan [11] contains a systematic
investigation of the high order differentiability of the norm function
in the Lebesgue-Bochner function spaces $L_p(X)$. Their results 
relate the continuously $n$-times differentiability of the norm on $X$ 
and of the norm on $L_p(X)$.

We improve the results in [11] in the following sense:
if $1<r\le p$ and the space $X$ admits
a norm  which is $C^n$-smooth, where $n$ is the
largest integer strictly less than $r$
and such that its  $n$-th  derivative is
$(r-n)$-H\"older on the unit sphere, then the norm on the space
$L_p(X)$ has the same properties of differentiability.

In the case $X=L_q$, the order of differentiability of the norm
given by the above result is, actually, the best possible
for bump functions and, consequently, also the best possible for
equivalent renormings.
For the results of non existence of high order differentiable bump
functions we investigate the situation in the subspace
$\bigoplus_{\ell_p} \ell_q^n$. As a consequence of our results we obtain 
that there is no twice Fr\'echet differentiable bump on the space
$\bigoplus_{\ell_2} \ell_4^n$, although this space has even a
twice G\^ateaux differentiable norm, as it was proved in [12].

This paper is divided into three sections. In Section I we introduce
functions with Taylor expansions of order $p$,
which are functions that verify Taylor's Formula around each point.
This notion is strictly weaker than the notion of 
$p$-times  differentiability in the case of $p$ being
an integer greater or equal
than $2$, and than
other categories of differentiability given in [2]. The main result in this
section is
Theorem 1-1, where we prove that if a Banach
space has modulus of convexity of power type $p$ and a bump function
with Taylor expansion of order $p$ at every point, then it has
a separating polynomial. This result improves Theorem 1 of [4],
and it relies on a different and much simpler proof,
using the variational principle of Stegall [16]. Applying this
to the classical space $L_p$ we improve the result of [2]
(Corollary 1-2). Some
variants are given in Theorem 1-5 for the case of spaces with
modulus of convexity with directional power type $p$.

Section II is devoted to results of high order differentiability of the
norm on $L_p(X)$. In the main result of this section, Theorem 2-1, 
starting with a $p$-homogeneous, $C^n$-smooth function on $X$ with
$n$-th derivative $(r-n)$-H\"older on the unit sphere
(where $1<r \leq p$ and $n$ is the largest
integer strictly less than  $r$) we construct 
a function on $L_p(X)$ with the same properties of differentiability.
In Theorem 2-3 we use this to prove the
result on higher order differentiability of the norm on $L_p(X)$
mentioned above. In Theorem 2-5 an analogue 
for bump functions is obtained.

Finally, in Section III, we give the best order of differentiability
of equivalent renormings on the space $L_p(L_q)$ in Theorem 3-1.
In order to do it,  we show in Theorem 3-2 that if the
space $\bigoplus_{\ell_p} \ell_q^n$ has a bump function with Taylor
expansion of order $p$ at every point, then $p,q$ are both even integers 
and $p$ is a multiple of $q$. In the proof of this Theorem, we obtain that
$\bigoplus_{\ell_p} \ell_q^n$ admits a separating polynomial
if and only if $p,q$ are both even integers and $p$ is a multiple of $q$.
We conclude from this that if the space $L_p(L_q)$ contains
$\bigoplus_{\ell_p} \ell_q^n$, (and this occurs except in trivial situations)
then $L_p(L_q)$ has a bump function with Taylor
expansion of
order $p$ at every point if and only if $p,q$ are both even integers
and $p$ is a multiple of $q$.

\bigskip
\noindent
{\bf I. Smooth norms in super-reflexive spaces}
\bigskip

Let $X$ be a real Banach space, let $p \geq 1$ be a real number
and let $f$ be  a real valued
function defined on $X$. We say that $f$ has a {\it Taylor
expansion of order} $p$ at the point $x\in X$, if there is
a polynomial $P$ of degree at most $n$, where $n=[p]$ is the integer part
of $p$, verifying
that
$$
\vert f(x+h) - f(x) - P(h) \vert = o(\|h\|^{p}).
$$
We say that $f$ is $T^p$-{\it smooth} if it
has a Taylor expansion of order $p$ at every point.
Note that if $f$ is $m$-times Fr\'echet differentiable on $X$,
then from Taylor's theorem we have that $f$ is $T^p$-smooth
for $1\leq p \leq m$.
In the same way as in [3], it is possible to
obtain some standard properties of $T^p$-smooth functions: in particular,
the composition of two $T^p$-smooth functions is again $T^p$-smooth.

\smallskip\noindent
As usual, a real-valued function on $X$ is said to be a {\it bump function} 
if it has bounded nonempty support.
We say that a polynomial $P$ on $X$ is a {\it separating polynomial}
if $P(0)=0$ and $P(x)\ge 1$ for all $x$ in the unit sphere of $X$.
It is known (and very easy to prove) that if $X$ admits a separating
polynomial,
then $X$ admits a $C^\infty$-smooth bump function.

\smallskip\noindent
Recall that the norm $\Vert.\Vert$ on a Banach space $X$ has {\it modulus
of convexity of power type} $p$ if there exists a constant $C>0$
such that for each $\varepsilon\in [0,2]$,
$$
\delta(\varepsilon) := \inf\{ 1-\Vert \frac{x+y}{2}\Vert: x,y\in X;\Vert
x\Vert\le 1;
\Vert y\Vert\le 1; \Vert x-y\Vert\ge\varepsilon\}\ge C\varepsilon^p .
$$

\bn{\bf Theorem 1.1} {\it Let $p\ge 1$ be a real number and let
$X$ be a Banach space. Let us assume :

1) There exists on $X$ a $T^p$-smooth bump function.

2) The norm of $X$ is uniformly convex with
modulus of convexity of power type $p$.

\noindent
Then there exists on $X$ a separating polynomial of degree $\le [p]$.}

\bigskip\noindent

Let $p\ge 1$ be a real number which is not an even integer,
and suppose that (for some measure space)
the space $L_p$ is infinite dimensional.
Kurzweil [10] proved that
$L_p$ does not admit any $C^{[p]}$-smooth bump function.
Bonic and Frampton [2] showed that in fact $L_p$ does not admit
any $B^p$-smooth bump function. Recall that,
according to Bonic and Frampton, a function $f$ is
said to be $B^p$-{\it smooth} if the
following is satisfied :

$(a)$ When $p$ is an integer, $f$ is $B^p$-smooth if it is $p$-times
Fr\'echet differentiable.

$(b)$ When $p$ is not an integer, $f$ is $B^p$-smooth if it is
$n$-times Fr\'echet differentiable, where $n=[p]$,
and for every $x$ in $X$
$$
\vert f^{(n)}(x+h)-f^{(n)}(x) \vert = o(\|h\|^{p-n}).
$$

\noindent
A real valued function on $X$ is said to be $H^p$-$smooth$ if
$f$ is $C^n$ on $X$, where $n$ is the largest integer strictly less
than $p$, and the $n$-th derivative $f^{(n)}$ is locally uniformly 
$(p-n)$-H\"older, that is, for
every $x \in X$, there are $\delta>0$ and $M>0$ such that
$$
\Vert f^{(n)}(y)-f^{(n)}(z)\Vert \leq M\Vert y-z\Vert^{p-n}
$$
for $y,z\in B(x;\delta)$. If the above inequality holds for
every $y,z \in X$ (with the same constant $M>0$), 
we say that $f$ is {\it uniformly} $H^p$-{\it smooth}. In the same
way, a $norm$ is said to be uniformly $H^p$-smooth if it is $C^n$
on $X -\{0\}$ and the $n$-th derivative is $(p-n)$-H\"older on the unit
sphere.

\noindent
It is clear that a uniformly $H^p$-smooth function need not be $T^p$-smooth,
but every
$B^p$-smooth function (and therefore every
$H^q$-smooth function for $q>p$)
is $T^p$-smooth; on the other hand, elementary examples on the real line show
that the converse is not true if $p\ge 2$. Therefore,
the following corollary improves the result of Bonic and Frampton.

\bigskip\noindent
{\bf Corollary 1.2. } {\it Let $p\ge 1$ be a real number which is not an
even integer.
Assume that $L_p$ is infinite dimensional. Then there is no
$T^p$-smooth bump function on $L_p$.}

\bigskip\noindent
{\bf Proof of Corollary 1.2:} Assume there exists a $T^p$-smooth bump on
$L_p$.
Since the modulus of convexity of the norm in $L_p$ is of power type $p$,
Theorem 1-1 shows that there is a separating polynomial on $L_p$, hence
there exists
a $C^\infty$-smooth bump on $L_p$, which contradicts the results
of Kurzweil, Bonic and Frampton given above.

\bigskip\noindent
For the proof of Theorem 1.1,
we shall use the following elementary lemmas from M. Fabian et al.

\bigskip\noindent
{\bf Lemma 1.3. [6] }{\it Let  $\delta(\varepsilon)$ be the modulus of
convexity of the norm
$\Vert.\Vert$ of $X$. Let $x, h\in X$ and $f\in X^*$ such that
$f(x)=\Vert x\Vert=\Vert f\Vert=1$, $f(h) = 0$ and $\varepsilon\le\Vert
h\Vert\le 2$.
Then :
$$
\Vert x+h\Vert\ge 1+\delta(\frac{\varepsilon}{2})
$$
Consequently, if the norm has modulus of convexity of power type $p$,
if $x, h\in X$ and $f\in X^*$ are such that
$f(x)=\Vert x\Vert\ne 0, \Vert f\Vert=1, f(h) = 0$ 
and $\Vert h\Vert\le 2\Vert x\Vert$, then :
$$
\Vert x+h\Vert - \Vert x\Vert \ge C\Vert x\Vert^{1-p}\Vert h\Vert^p
$$
}

\bigskip\noindent
{\bf Lemma 1.4. [6] }{\it Let $k\ge 2$ be an integer and let $F$ be a
finite codimensional
subspace of a Banach space $X$. Assume that $P$ is a polynomial on $X$ of degree
$\le k$ which is a separating polynomial on $F$.
Then there is a separating polynomial of degree $\le k$ on $X$.}

\bigskip\noindent
{\bf Proof of Theorem 1.1:}
Let $b$ be a bump function on $X$ with
Taylor expansion of order $p$ at any point and
such that $b(0)=0$.
Let  $ \varphi : X \to \reel \cup\{+\infty\}$ be the function
defined by :
$\displaystyle\varphi(x)=\frac{1}{b(x)^2}$ if $b(x)\ne 0$
and $\varphi(x)=+\infty$ otherwise. The function
$\varphi-\Vert .\Vert$ is lower
semicontinuous,
bounded below and identically equal to $+\infty$ outside a bounded set.
On the other hand, since $X$ is uniformly convex, it has the Radon-Nikodym
property.
According to Stegall variational principle [16],
there exists $g\in X^*$ (actually
a dense $\Cal G_\delta$ in $X^*$ of such $g$'s) such that $\varphi-\Vert
.\Vert + g$
attains its minimum at some point $x$. So for every $h\in X$ :
$$
\varphi(x+h) - \Vert x+h\Vert + g(x+h)\ge \varphi(x) - \Vert x\Vert + g(x)
$$
Consequently :
$$
\varphi(x+h) - \varphi(x) + g(h) \ge \Vert x+h\Vert - \Vert x\Vert
$$
Next $x\ne 0$ because $\varphi(0)=+\infty$.
Let $f\in X^*$ such that $\Vert f\Vert = 1$ and $f(x)=\Vert x\Vert$.
Using lemma 1-3, for every $h\in Ker(f)\cap Ker(g)$ :
$$
\varphi(x+h) - \varphi(x) \ge \Vert x+h\Vert - \Vert x\Vert
\ge C\Vert x\Vert^{1-p}\Vert h\Vert^p
$$
If we denote $C(x) = C\Vert x\Vert^{1-p}$, we thus have :
$$
\varphi(x+h) - \varphi(x) \ge C(x)\Vert h\Vert^p
$$
We now use the fact that $\varphi$ has
a Taylor expansion of order $p$ at $x$ :
there exists a polynomial $P$ of degree $\le [p]$
and a function $R$
such that $\varphi(x+h)-\varphi(x) = P(h) + R(h)$, $P(0) = 0$ and
$\displaystyle\lim_{h\to 0} \frac{R(h)}{\Vert h\Vert^p} = 0$. We fix
$\varepsilon > 0$ such that
$\displaystyle\vert R(h)\vert\le\frac{C(x)}{2}\Vert h\Vert^p$
whenever $h\in X, \Vert h\Vert = \varepsilon$. Therefore,
if $h\in Ker(f)\cap Ker(g)$ and $\Vert h\Vert = \varepsilon$, then :
$$
P(h)\ge\frac{C(x)}{2}\Vert h\Vert^p
$$
The polynomial $Q$ defined by
$$
Q(h) = \frac{2}{C(x)\varepsilon^p}  P(\varepsilon h)
$$
is a separating polynomial of degree $\le [p]$ on $Ker(f)\cap Ker(g)$.
By lemma 1-4, there is a separating polynomial of degree $\le [p]$ on $X$.

\bigskip\noindent
Variants of  Theorem 1-1 are possible : Recall that
the norm $\Vert.\Vert$ on a Banach space $X$ has modulus
of convexity of directional power type $p$ if for every finite dimensional
subspace $F\subset X$, there exists a constant $C_F>0$
such that for each $\varepsilon\in[0,2]$,
$$
\delta_F(\varepsilon) := \inf\{ 1-\Vert \frac{x+y}{2}\Vert; x,y\in F;\Vert
x\Vert\le 1;
\Vert y\Vert\le 1; \Vert x-y\Vert\ge\varepsilon\}\ge C_F\varepsilon^p .
$$
The above proof shows the following result :

\bn{\bf Theorem 1.5.} {\it Let $p\ge 1$ be a real number and let
$X$ be a Banach space. Let us assume :

\smallskip
(1) There exists on $X$ a
bump function with Taylor expansion of order $p$ at any point.

(2) The norm of $X$ is uniformly convex with
modulus of convexity of directional power

\quad\,\,\,\,
type $p$.

\smallskip\noindent
Then there exists on $X$ a polynomial $Q$ (of degree $\le p$) such that
$Q(h)>0$
whenever $h\in X$ and $\Vert h\Vert = 1$.}

\bigskip\noindent
{\bf Proof:} Indeed, following the proof of Theorem 1-1, We obtain $f,
g\in X^*$,
$\varepsilon>0$, a constant $C(x,F)>0$
and a polynomial $P$ of degree $[p]$ such that $P(0)=0$ and
$P(h)\ge C(x,F)\Vert h\Vert^p$ whenever $h\in F\subset Ker(f)\cap Ker(g)$ and
$\Vert h\Vert = \varepsilon$. The polynomial $Q$ defined by
$Q(h)=P(\varepsilon h)
+f^2(h)+g^2(h)$ satisfies $Q(0)=0$ and $Q(h)>0$ whenever  $h\in X$ and
$\Vert h\Vert = 1$.

\bn
The following variant of Theorem 1-1 will be used in section 3.

\bn
{\bf Theorem 1.6. }{\it Let $1<p<\infty$, and
let $X$ be a Banach space with the Radon-Nikodym property
satisfying the following property :
\smallskip
``For all $x\in X$ there is a constant $C>0$, such that
for each $\delta >0$ there is a finite codimensional
subspace $H_{\delta}$ of $X$, such that
$\Vert x + h \Vert - \Vert x\Vert \geq C \Vert h \Vert^{p}$
for all $h\in H_{\delta} $ with $\Vert h \Vert = \delta $."

If $X$ has a $T^p$-smooth bump function, then it has a separating polynomial.}

\bn
{\bf Proof : }
We proceed as in Theorem  1-1 :
Let $b$ a bump function which is $T^p$-smooth
and such that $b \geq 0$ and $b(0)=0$.
We consider
$\varphi$  defined in the proof of Theorem 1-1, and, proceeding in the same
way as there,
there exists $g\in X^{*}$ such that $\varphi - \Vert . \Vert + g$
attains its minimun at some point $x \in dom(\varphi)$.
Therefore, for every $h \in  Ker(g)$ :
$$
\varphi(x+h) - \varphi(x)
\geq \Vert x+h \Vert - \Vert x \Vert.
$$
For such a point $x$, there is a constant $C>0$ verifying
that for each
$\delta >0$ there is a finite codimensional
subspace $H_{\delta}$ of $X$ such that
$$
\varphi(x+h) - \varphi(x) \geq C\Vert h \Vert^{p},
$$
for all $h \in H_{\delta} \cap Ker(g)$ with $\Vert h \Vert =\delta$.
Next,  we use that the function $\varphi$ has a Taylor expansion of
order  $ p$ at $x$. Then we may choose  some $0< \delta \leq 1 $ such that
if $h \in X$ with $\Vert h \Vert \leq \delta$ :
$$
\vert
\varphi(x+h) - \varphi(x) - P(h)
\vert \leq \frac{C}{2} \Vert h \Vert^{p}
$$
where $P$ is a polynomial of degree $\le [p]$ vanishing at the origin.
Therefore, for $h \in H_0 = H_{\delta} \cap Ker(g)$, with
$\Vert h \Vert =\delta$ we have that
$$
\vert P(h) \vert \geq
\varphi(x+h) - \varphi(x) - \frac{C}{2} \delta^p
\geq \frac{C}{2} \delta^p.
$$
From the above the existence of a
separating polynomial on $H_0$ follows. Since $H_0$ is
a finite codimensional
subspace of $X$, by Lemma 1-4
the space $X$ has a separating polynomial, as we required.

\bigskip
\noindent
{\bf II. Smoothness of $L_p(X)$}
\bigskip
In what follows, let $1<p< \infty$ and let
$(\Omega, \mu)$ a measure space such that the corresponding $L_p$-space
is infinite-dimensional. Recall that for a Banach space $X$
the function space $L_p(X)$ is defined as follows :
$$
L_p(X) = \big \{ u : \Omega \to X / u \text{ measurable and } 
\Vert u \Vert_p = 
\left( \int_{\Omega} \Vert u(s)) \Vert^p d \mu (s) \right)^{1/p}
< \infty \big \}.
$$

\bn
As usual, a real valued function $f$ on $X$
is said to be $p$-homogeneous if $f(\lambda x)=\vert \lambda \vert^p f(x)$
for every $x \in X$ and every $\lambda \in \reel$.
The main result of this section is the following:

\bn
{\bf  Theorem 2.1. }
{\it Let $X$ be a Banach space, let $1<r \leq p$ and
consider $n$ the largest integer strictly less than $r$. Let $f$ be
a $p$-homogeneous real function defined on $X$. Suppose
that  $f$ is $C^n$
on $X -\{0\}$ and $f^{(n)}$ is $(r-n)$-H\"older on the unit sphere of $X$.
Then, the function
$\widehat{f}$ on $L_{p}(X)$ defined by
$$
\widehat{f}(u) = \int_{\Omega} f(u(s))d \mu(s)
$$
is a $p$-homogeneous $C^{n}$-smooth function on $L_p(X)$ and
$\widehat{f}^{(n)}$ is $(r-n)$-H\"older on bounded subsets of $L_p(X)$. }

\bigskip\noindent
We first establish the following lemma :

\bn
{\bf Lemma 2.2.}
{\it Let $X$ be a Banach space, let $1<r \leq p$ and
consider $n$ the largest integer strictly less than $r$. Let $f$ be
a $p$-homogeneous real function defined on $X$
and suppose that  $f$ is  $C^n$
on $X -\{0\}$ and $f^{(n)}$ is $(r-n)$-H\"older on the unit sphere of $X$.

Then $f$ is $C^n$-smooth on $X$, $f^{(j)}(0)=0$ for $j=1,2, \dots ,n$ and
there exist constants $C_1, C_2, \dots , C_n >0$ so that :

\medskip
(1) $\Vert f^{(n)}(x) - f^{(n)}(y) \Vert \leq
C_n (max \{\Vert x \Vert, \Vert y \Vert \})^{p-r} \Vert x-y \Vert^{r-n}$
for every $x, y \in X$.

\medskip
(2) $\Vert f^{(j)}(x) - f^{(j)}(y) \Vert \leq
C_j (max \{\Vert x \Vert, \Vert y \Vert \})^{p-(j+1)} \Vert x-y \Vert$
for every $x, y \in X$ and $j=1,2, \dots , n-1$.
}

\bn
{\bf Proof of lemma 2.2 : }
Since $f$ is $p$-homogeneous, its derivatives are positively
homogeneous, that is, $f^{(j)}(\lambda x) = \lambda^{p-j}
f^{(j)}(x)$ for every $\lambda >0, x \neq 0$ and $j=1,2,\dots,n$.

Since $f^{(n)}$ is $(r-n)$-H\" older on the unit sphere of $X$, there exists
a constant $K_n>0$ such that
$$
\Vert f^{(n)}(x) - f^{(n)}(y) \Vert \leq K_n \Vert x-y \Vert^{r-n}
$$
whenever  $\Vert x \Vert = \Vert y \Vert =1$. Now, fix $x,y$ in $X$
with $0 < \Vert y \Vert \leq
\Vert x \Vert$ ,  and consider $z= \frac{\| y \|}
{\| x \|} x $. Then, since $\Vert z \Vert = \Vert y \Vert$,
$$
\Vert f^{(n)}(z) - f^{(n)}(y)\Vert =
\big\Vert f^{(n)} \bigl(\frac{z}{\Vert z \Vert}\bigr) -
f^{(n)} \bigl(\frac{y}{\Vert y \Vert}\bigr) \big\Vert \Vert  y \Vert^{p-n}
\leq K_n \Vert z-y\Vert^{r-n} \Vert y \Vert^{p-r};
$$
and
$$
\Vert z-y \Vert = \frac{1}{\Vert x \Vert} \Big\Vert x \|y \| - \|x \| y
\Big\Vert \leq  \Big\vert \| y \| - \| x \| \Big\vert + \|x-y \|
\leq 2\|x-y\| .
$$
Hence,
$$
\| f^{(n)} (z) - f^{(n)} (y)\|
 \leq 2^{r-n} K_n \|y \|^{p-r} \|x-y\|^{r-n} \leq
$$
$$
\leq K'_n (max\{ \|x\|,\|y\|\})^{p-r} \|x-y\|^{r-n} \qquad (*)
$$
Now,
$$
\| f^{(n)} (x) - f^{(n)} (z) \| = \Bigl(1- \big ( \frac{\|y \|}{\|x\|} \big
)\Bigr)^{p-n}
\| f^{(n)} (x) \| = (\|x\|^{p-n}  -\|y\|^{p-n} )
\| f^{(n)} \bigl(\frac{x}{\|x\|}\bigr) \|.
$$
The function $f^{(n)}$ is $(r-n)$-H\"older on $S_{X}$,  and
consequently it is uniformly bounded on the unit sphere; hence,
by homogeneity, there is a constant $D_{n}>0$ such that
$\| f^{(n)}  (x) \| \leq D_n \|x \|^{p-n}$ for all $x\neq 0$.
\smallskip
Now, we consider two cases : either $p-n>1$,  or $0<p-n\leq 1$.
In the first case,
By the mean value theorem applied to the real function
$\varphi(t) = t^{p-n}$, we have
$$
\Big\vert \|x\|^{p-n} - \|y\|^{p-n} \Big\vert \leq (p-n)
(max\{\|x\|,\|y\|\})^{p-n-1}
\|x-y\| \leq
$$
\smallskip
$$
\leq (p-n) (max\{\|x\|,\|y\|\})^{p-n-1}(\|x\|+\|y\|)^{1-r+n} \|x-y\|^{r-n}
$$
\smallskip
$$
\leq 2^{1-r+n} \max\{\|x\|,\|y\|\}^{p-r}
\|x- y \|^{r-n}.
$$
Otherwise, if $0<p-n\leq1$, then
$$
\Big\vert  \|x\|^{p-n} -\|y\|^{p-n} \Big\vert \leq
\|x-y\|^{p-n} \leq
(\|x\|+\|y\|)^{p-r}\|x-y\|^{r-n}
$$
\smallskip
$$
\leq 2^{p-r} (\max\{\|x\|,\|y\|\})^{p-r}\|x-y\|^{r-n}.
$$
Hence, in both cases we have that for some constant
$K_{n}''$,
$$
\| f^{(n)}(x) - f^{(n)}(z)\| \leq K_{n}''
(max \{\|x\| ,\|y\|\})^{p-r} \| x-y \|^{r-n}    (**)
$$
\bigskip
From (*) and (**) we obtain (1). As a consequence of it,
we obtain (2) for all $x,y \neq 0$.  Indeed, consider
first $j=n-1$ and $x,y\neq 0$. In the case that
$0 \notin [x,y]$, by the mean value theorem applied to the
function $f^{(n-1)}$ we have :
$$
\| f^{(n-1)}(x) - f^{(n-1)}(y) \| \leq \sup_{ w \in [x,y]}
\| f^{(n)}(w) \|. \|x-y\| \leq
$$
\smallskip
$$
\leq D_n \max\{ \|x\|,\|y\|\}^{p-n} \|x-y\|.
$$
\smallskip
Otherwise, if $0\in [x,y]$, then
the linear span of $x$ coincides whith that of $y$; assuming that
$dim X \geq 2$, we may choose $z\in X$ such that $z$ does not
belong to the linear span of $x$ and verifies that
$\|z\| = \min \{\|x\|,\|y\|\} >0$. Then,
$$
\| f^{(n-1)}(x) - f^{(n-1)}(y)\| \leq
\| f^{(n-1)}(x) - f^{(n-1)}(z)\| + \| f^{(n-1)}(z) - f^{(n-1)}(y)\| \leq
$$
and since $0\notin [x,z]\cup [y,z]$, applying again the mean value theorem,
$$
\leq D_{n} \max\{ \|x\|,\|y\|\} ^{p-n}(\|x-y\|+\|z-y\|)
$$
and since $\|x-z\|,\|z-y\| \leq \|x-y\|$, we have $(2)$ for $j=n-1$.
In an analogous way, using that $f^{(n-1)}$ is bounded in $S_{X}$
(it is indeed Lipschitz), and proceeding as above, $(ii)$ is obtained for
$x,y\neq 0$.
\bigskip
In particular, each $f^{(j)}$ is bounded in $S_{X}$; and by homogeneity,
for all $x\neq 0$,
$$
\| f^{(j)}(x) \| \leq D_j \|x\|^{p-j}.
$$
This implies that $f$ is $n$-times differentiable at $0$ and
that $f^{(j)}(0)=0$ for $j=1,\dots,n$. Therefore $(1)$ and $(2)$ hold
indeed for all $x,y \in X$, and consequently $f\in C^{n}(X)$.
Finally, the result is clear if $dim(X) =1$, since all
$p$-homogeneous function on $\reel$ are of the form $f(t) = a\vert t\vert^{p}$
for some $a \in \reel$.

\bn
{\bf Proof of  Theorem 2.1:}
The function taking $s \to f(u(s))$ is measurable, since $f$ is a continuous
function. On the other hand, since $\vert f(x)\vert \leq C \|x\|^{p}$
for all $x\in X$, the function $\widehat{f}$ is well defined; indeed,
$$
\int_{\Omega} \vert f(u(s))\vert d\mu (s) \leq
\int_{\Omega} C \| u(s)\|^{p} d\mu (s) < \infty.
$$
Now, for each $u \in L_{p}(X)$ and each $j=1,\dots,n$, let
$P_j(u)$ be the
$j$-homogeneous polynomial on $L_p(X)$ defined by
$$
P_j(u)(h) = \int_{\Omega}  f^{(j)}\bigl(u(s)\bigr)\bigl(h(s)\bigr)  d \mu (s)
$$
for all $h \in L_p(X)$. These
polynomials are well defined; indeed, $s \to f^{(j)}(u(s))(h(s))$
is a measurable function and, since $ \| f^{(j)} (x) \| \leq C_j
\|x\|^{p-j}$  for
all $x\in X$, we have
$$
\int_{\Omega} \big\vert f^{(j)}\bigl(u(s)\bigr)\bigl(h(s)\bigr) \big\vert d
\mu (s) \leq
\int_{\Omega} C_{j}\| u(s) \|^{p-j} \|h(s)\|^{j} \leq
$$
and by H\"older's inequality
$$
\leq C_j (\int_{\Omega} \| u(s)\|^p d  \mu(s) )^{\frac{p-j}{p}}
(\int_{\Omega} \|h(s)\|^{p})^{\frac{1}{p}} < \infty .
$$
\bigskip
In order to obtain that $\widehat{f}$ is $C^{(n)}$-smooth on
$L_p(X)$ and $\widehat{f}^{(j)} = P_j$, we apply the
converse of Taylor's Theorem (see e.g. [1, Th. 4-11]). To apply this,
we need to prove that $\widehat{f}$ has a Taylor
expansion of order $n$ at any $u\in L_p(X)$ (with
the polynomial $\sum_{j=1}^{n} \frac{1}{j!}P_j(u)$) and,
also,
that the function $u \to P_{j}(u)$ is continuous for all
$j=1,\dots,n$. We begin with the second requirement. Consider
$1\leq j \leq n$ and $u,v \in L_p(X)$; then :
$$
\vert (P_j(u) - P_j(v))(h) \vert \leq
\int_{\Omega} \|f^{(j)} (u(s))-f^{(j)} (v(s))\| \|h(s)\|^{j} d \mu (s) \leq
$$
and by (ii) of Lemma 2-2,
$$
\leq C_j\int_{\Omega} max\{ \|u(s)\|,\|v(s)\|\}^{p-j-1} \|u(s) - v(s)\| 
\|h(s)\|^{j} d \mu (s) \leq
$$
and by H\"older's inequality, if we consider
the function
$\Psi \in L_p(\Omega)$ defined
by $\Psi(s) =max \{ \|u(s)\|,\|v(s)\|\} $, then
$$
\leq C_j \|\Psi\|^{p-j-1} \|u-v\| \|h\|^{j}
$$
$$
\leq C_j (\|u \|^p + \|v \|^p)^{\frac{p-j-1}{p}} \|u-v\| \|h\|^{j}.
$$

Therefore, the function $u \to P_j(u)$ is Lipschitz on bounded subsets
for $1\leq j\leq n-1$,
and in particular,  it is continous.  Analogously,
$$
\vert ( P_n(u) - P_n(v) ) (h) \vert \leq C_n \| \Psi \|^{p-r} \|u-v\|^{r-n}
\|h\|^{n}
$$
$$
\leq C_n ( \| u \|^p +\| v \|^p)^{\frac{p-r}{p}} \|u-v\|^{r-n} \|h\|^{n}, 
$$
which imply that $u \to P_{n}(u)$ is $(r-n)$-H\"older on
bounded subsets of $L_{p}(X)$.
\bigskip
To obtain a Taylor expansion for $\widehat{f}$, we apply 
Taylor's formula with integral remainder; hence,
for $x,h\in X$,
$$
\vert  f(x+h) - f(x) - \sum_{j=1}^{n} \frac{1}{j!}f^{(j)}(x)(h)
\vert =
$$
$$
= \vert \int_{0}^{1} \frac{(1-t)^{n-1}}{(n-1)!}(f^{(n)}(x+th) - 
f^{(n)}(x))(h) dt \vert 
\leq \frac{C_n}{n!} (\|x\| + \|h\|)^{p-r}\|h\|^{r} .
$$

Now denote
$$
R_n (u;h) = \widehat{f}(u+h) - \widehat{f}(u) - \sum_{j=1}^{n} P_j(u)(h) .
$$
Then
$$
\vert R_n (u;h) \vert 
=\big\vert \int_{\Omega} \Bigl( f(u(s) + h(s) ) - f(u(s) - \sum_{j=1}^{n}
\frac{1}{j!}f^{(j)}(u(s))(h(s)) \Bigr) d \mu (s) \big\vert
$$
$$
\le \frac{C_n}{n!} \int_{\Omega} (\|u(s)\| + \|h(s)\|)^{p-r}\|h(s)\|^{r} \le
\frac{C_n}{n!} \Bigl(\int_{\Omega} (\|u(s)\| + \|h(s)\|)^{p} \Bigr)^\frac{p-r}{p}
\Bigl(\int_{\Omega} \|h(s)\|^{p}\Bigr)^{\frac{r}{p}}
$$
The last inequality above follows from H\"older's inequality.
Since $s \to \|u(s)\| + \| h(s)\|$ belongs to $ L_p(\Omega) $,
$$
\vert R_n (u;h) \vert 
\leq \frac{C_n}{n!} (\| u \| + \| h \|)^{p-r} \| h \|^{r}.
$$
Therefore we have that, for each $u_0 \in L_p(X)$, 
$$
\lim_{(u,h) \to (u_0,0)} \frac{R_n(u;h)}{\| h \|^n} =0,
$$
and the result follows.

\bn
Now we apply the above Theorem to the particular case $X= L_q$
(for a probably different measure space). We obtain the
following  results concerning the high order differentiability of the norm on
$L_p(L_q)$.

\bn
{\bf Corollary  2.3. }{\it Let $1<p,q< \infty$  and $X=L_p(L_q)$.
Consider $r=min\{p,q\}$. Then :

If $r$ is not an even integer, or $r=p$ is an even integer,
then the norm on $X$ is uniformly $H^r$-smooth.

If $r$ is an even integer and $r=q$ then the norm in $X$ is
uniformly $H^p$-smooth.
Moreover, if $p$ is a multiple of $q$ then the norm is even
$C^{\infty}$-smooth. }

\medskip\noindent
{\bf Proof of Corollary 2.3:}
Since $r \leq q$, the norm on $L_q$ is  uniformly $H^r$-smooth [2]; by
Theorem 2-1 the function
$$
F(u) = \int_{\Omega}\big( (\Vert u(s) \Vert_q)^p \big)d \mu(s) 
$$
is $p$-homogeneous, $C^n$-smooth and such that $F^{(n)}$ is
$(r-n)$-H\"older on $S_{X}$, where
$n$ is the largest integer strictly less than $r$. Therefore, the
norm on $X$ is uniformly $H^r$-smooth.

If $r=q$ is an even integer, the norm in $L_r$ is
$C^{\infty}$-smooth, and in
particular uniformly $H^p$-smooth. Using again Theorem 2-1, the norm
in $X$ is uniformly $H^p$-smooth. Moreover, if $p$ is a multiple of $q$
then the expression
$F$ defines a $p$-homogeneous separating polynomial on $X$
and the norm is, in particular, $C^{\infty}$-smooth.
\bigskip
From Theorem 2-1 we have that if the norm of $X$ is uniformly 
$H^r$-smooth, and $1<r\leq p$, then the norm of $L_p(X)$ is also 
uniformly $H^r$-smooth. Next we obtain the analogous result for 
bump functions. To this end, we use the following construction given 
in [5], which produces
an homogeneous smooth function from a smooth bump function:

\medskip\noindent
{\bf Proposition  2.4. }{\it Let $1<r\leq p$ and let $n$ be the largest
integer strictly less than $r$. Assume that a Banach space $X$ admits an
uniformly
$H^r$-smooth bump function. Then, there is a $p$-homogeneous
function $f$ on $X$ such that $f$ is $C^n$-smooth on $X$,
$f^{(n)}$ is $(r-n)$-H\"older on bounded subsets,
and there exist some constants $A,B>0$, such that
$$
A \Vert x \Vert^p   \leq f(x) \leq B \Vert x\Vert^p
$$
for all $x\in X$.}

\bn
{\bf Proof : } Let $\varphi$ be a uniformly $H^r$-smooth bump function on $X$.
We proceed as in [5, Proposition II-5-1] : Consider $\widehat{\varphi}
=1-exp(-\varphi^2)$ and
$\widehat{\psi}$ defined by
$$
\widehat{\psi}(x) = \int_{-\infty}^{\infty}\widehat{\varphi}(sx)ds.
$$
Set $\psi(x) = \frac {1}{\widehat{\psi}(x)}$ if $x\neq 0$ and
$\psi(0)=0$. Then, the function $\psi$ is $1$-homogeneous, $C^1$-smooth
on $X-\{ 0 \}$, and 
$$
a\Vert x \Vert   \leq \psi (x) \leq  b\Vert x\Vert;
$$
where $a,b$ have been chosen verifying that
$\widehat{\varphi}(x) > \frac{1}{2} \widehat{\varphi}(0)$ and
$\widehat{\varphi}(x)=0$ if $\Vert x \Vert \geq b$.

Now,  given $x-0 \neq 0$, there are $\varepsilon >0$ and  $K>0$, such that
if $\Vert x-x_0 \Vert< \varepsilon$, then
$$
\widehat{\psi}(x) =  \int_{-K}^{K}\widehat{\varphi}(sx)ds.
$$
Then, define for $x \in B(x_0;\varepsilon)$, $h \in X$ and $j=1,2,\dots,n$:
$$
P_j(x)(h) = \int_{-K}^{K} s^j \widehat{\varphi}^{(j)}(sx)(h)ds.
$$
Each mapping
$x \to P_j(x)$ is continuous on $B(x_0;\varepsilon)$ and
$$
\vert \widehat{\psi}(x+h) - \widehat{\psi}(x) - \sum_{j=1}^{n} \frac{1}{j!}
P_j(x)(h) \vert = o(\Vert h \Vert^n),
$$
uniformly for $x \in B(x_0;\varepsilon)$.
Then, by using the converse of Taylor's theorem (see [1, Th. 4-11]) 
we obtain that $\widehat \psi$ is
$C^n$-smooth and $\widehat \psi^{(j)} = P_j$ on $B(x_0;\varepsilon)$ for 
$j=1,\dots,n$.
Therefore, $\widehat \psi$ is
$C^n$-smooth on $X-\{0\}$.
In order to prove that $\widehat{\psi}^{(n)}$ is  $(r-n)$-H\"older
on the unit sphere note that, given $x \in S_X$, if
$\vert s \vert \geq b$ then
$\widehat{\varphi}^{(n)}(sx)=0$; therefore,
$$
\widehat{\psi}^{(n)}(x)(h) =  \int_{-b}^{b} s^n
\widehat{\varphi}^{(n)}(sx)(h)ds .
$$
Hence, for $x,y \in S_X$ and $h \in X$:
$$
\widehat{\psi}^{(n)}(x)(h) -
\widehat{\psi}^{(n)}(y)(h)  =
\int_{-b}^{b} s^n \big( \widehat{\varphi}^{(n)}(sx) 
-\widehat{\varphi}^{(n)}(sy)
\big)(h)ds.
$$
Using that $\widehat{\varphi}$ is uniformly $H^r$-smooth, 
we have that there is a constant $M>0$, such that
$$
\Vert \widehat{\psi}^{(n)}(x) -
\widehat{\psi}^{(n)}(y) \Vert \leq M \Vert x-y \Vert^{r-n}.
$$
Therefore, the $n$-th derivative of $\widehat{\psi}$ is $(r-n)$-H\"older on
the unit sphere, and so is the $n$-th derivative of $\psi$. Consider now
$f(x) = \psi(x)^p$. By lemma 2.2, the function $f$ satisfies the required
conditions.

\bn
{\bf Theorem  2.5. }{\it Let $1<r \leq p < \infty$ and let $X$ be a Banach
space. Then,
the following statements are equivalent:

1) $X$ admits a uniformly $H^r$-smooth bump function.

2) $L_p(X)$ admits a uniformly $H^r$-smooth bump function.

3) $L_p(X)$ admits an $H^r$-smooth bump function.}

\bn
{\bf Proof: }
First we will prove that (1) implies (2). Assume that $X$
admits a uniformly $H^r$-smooth bump function,  
and let $n$ be the largest integer strictly less than $r$.
According to
Proposition 2-4 there is a $C^n$-smooth, $p$-homogeneous function $f$ on $X$,
such that $f^{(n)}$ is
$(r-n)$-H\"older on bounded subsets, and for some constants
$A,B >0$ we have
$$
A \Vert x \Vert^p   \leq f (x) \leq B \Vert x\Vert^p.
$$
Consider the function on $L_p(X)$:
$$
\widehat{f}(u) = \int_{\Omega} f(u(s))d \mu(s)
$$
By Theorem 2-1, $\widehat{f}$ is $C^n$-smooth on $L_p(X)$
and $\widehat{f}^{(n)}$ is $(r-n)$-H\"older on bounded subsets.
By composing with a convenient $C^{\infty}$-smooth function
$\varphi: \reel \to \reel$, we obtain a
uniformly $H^r$-smooth bump function on $L_p(X)$.

Since (2) implies (3) is trivial, we now prove that (3) implies (1). 
Assume that $L_p(X)$ admits a $H^r$-smooth bump function. 
We claim that $X$ does not
contain an isomorphic copy of $c_0$. Indeed: otherwise, the subspace
$Y = \bigoplus_{\ell_p} c_0^n$ of $L_p(X)$
would have an $H^r$-smooth bump function, and
therefore a uniformly $H^r$-smooth bump function, since
it does not contain an isomorphic copy of
$c_0$ (see [5, Theorem V-3-1]); in particular,
$Y$ would be superreflexive, which is not possible.
Therefore, since $X$ does not contain an isomorphic copy of $c_0$, it is
uniformly $H^r$-smooth (see [5, Th. V-3-1]) as we required.

\bigskip \noindent
{\bf III. Best order of smoothness of $L_p (L_q)$ }

\bigskip

In this section we consider infinite-dimensional spaces $L_p$ and $L_q$
over probably different measure spaces.
We prove that the order of differentiability of the norm
in the space $L_p(L_q)$ given in Section II is, actually, the best possible
for bump functions, and, in particular, for equivalent renormings.

\bn
{\bf Theorem 3.1. }{\it Let $1<p,q< \infty$  and $X=L_p(L_q)$.
Consider $r=min\{p,q\}$. Then :

If $r$ is not an even integer, or $r=p$ is an even integer,
then there is no $T^r$-smooth bump function on $X$.

If $r$ is an even integer and $r=q$ then there is a $T^p$-smooth
bump function on $X$ if and only if only both $p,q$ are even integers and
$p$ is a  multiple of $q$.
Moreover, in this case, the space $X$ has a separating polynomial.}
\bigskip

In order to prove Theorem 3-1 we investigate the best order
of smoothness of a certain subspace of $L_p(L_q)$, namely the
subspace $\bigoplus\Sb \ell_{p} \endSb\ell_{q}^{n}$.

\bn
{\bf Theorem 3.2. }{\it
Let   $1<p,q<\infty $ and   $X = \bigoplus_{\ell_p} \ell_q^n$. 
Assume that  $X$ has a $T^p $-smooth bump function. Then $p,q$ are both 
even integers and $p$ is a multiple of $q$. 
Moreover, $X$ admits a separating polynomial.}

\bn
For the proof of this theorem we need the following lemma :

\bn
{\bf Lemma 3.3. } {\it Let $N \in \entier$ and
assume that $X$ is a separable Banach space such
that every weakly null and normalized sequence
in $X$ has a subsequence which is equivalent to
the unit vector basis in $\ell_{N}$. If
$X$ has a separating polynomial of
degree at most $N$, then there is an $N$-homogeneous
separating polynomial. }

\bn
{\bf Proof : } First note that $X$ is superreflexive and
therefore  $X^{*}$ is also separable.
Consider $\{ z_{n}^{*} \}$ a countable dense set on the unit sphere of
$X^{*}$. Now let $P$ be a separating polynomial
of degree $N$, with $P=P_{1} +\dots + P_{N-1}
+P_{N}$, where each $P_{i}$ is an $i$-homogeneous
polynomial.

Assume that there is no $N$-homogeneous separating
polynomial on $X$. Then,
$$
\inf \Sb \Vert x \Vert =1 \endSb P_{N}(x) = 0
$$
and this infimum is also $0$ on the unit sphere of
each finite codimensional subspace. Then, we can construct
a sequence in the following way : for each
$ n \in \entier$ consider $x_{n} \in Ker z_{1}^{*} \cap \dots
\cap Ker z_{n-1}^{*}$
with $\Vert x_n \Vert =1$ and verifying that
$\vert P_N(x_n)\vert \leq \frac{1}{n}. $
By the construction, the sequence $\{x_n\}$ is weakly null,
and then it has a subsequence
which is equivalent to the unit vector basis of $\ell_N$. In particular
it has a subsequence $\{x_{n_k}\}$
with an upper $N$-estimate, and by [8] for instance, it follows that
$\{(P_{1}+\dots +P_{N-1})(x_{n_k})\}_k \to 0$.
Hence, $\{P(x_{n_k})\}$ also converges to zero,
which contradicts the fact that $P$ is separating.
\bigskip
\noindent
{\bf Proof of Theorem 3.2:} The proof will be done in two steps:
\bn
{\bf Step 1. }{\it
Let   $1<p,q<\infty $
and   $X = \bigoplus_{\ell_p} \ell_q^n$. Assume that  $X$ has a
$T^p $-smooth bump function. Then $X$ admits a separating polynomial. }
\bn
{\bf Proof of Step 1:} It is enough to prove that the space $X$ verifies
the hypothesis
of Theorem 1-6, that is, the property:

``For all $x\in X$ there is a constant $C>0$, such that
for each $\delta >0$ there is a finite codimensional
subspace $H_{\delta}$ of $X$, such that
$\Vert x + h \Vert - \Vert x \Vert \geq C \Vert h \Vert^{p}$
for all $h\in H_{\delta} $ with $\Vert h \Vert = \delta $."

If $x=0$ the result is trivial. Now supose $x \neq 0$.
For each $n \in \entier$, we denote by
$\{e_1^{n},\dots, e_n^{n}\}$ the usual basis
in $\ell_q^n$. The point $x$ can be expressed
as $x= \{ x_n\}$, where
each $x_n = \sum_{i=1}^{n} x_i^{n} e_i^{n}$, and
$
\Vert x \Vert^{p} = \sum_{i=1}^{\infty} \Vert x_n \Vert_{\ell_n^q}^{p}.
$
We denote $x^{n} = x_1 + \dots + x_n$, and it is clear that
$\Vert x- x^n \Vert \to 0$.
For each $n$ we consider the following subspaces of $X$
$$
H^{n} = \bigoplus_{i=1}^{n} \ell_q^i \qquad\text{and}\qquad
H_n = \bigoplus_{i=n+1}^{\infty} \ell_q^{i}.
$$
It is easy to check that if $y \in H^n$ and $z \in H_n$ then,
$$
\Vert y + z \Vert^p = \Vert y \Vert^p + \Vert z \Vert^p.
$$
Now, applying the mean value theorem to the
real function $\lambda(t) = t^p$ on the interval
$I=[\Vert x \Vert/2, \Vert x \Vert +1]$, there is
a positive constant $C>0$, such that if $s,t \in I$ and
$s<t$, the following inequality holds :
$$
t - s \geq C (t^p - s^p).  \qquad (*)
$$
Therefore, if $h\in H_n$ and $\Vert h\Vert=\delta\leq 1$, we obtain
from (*) that
$$
\Vert x+h \Vert - \Vert x \Vert \geq
\Vert x^n + h \Vert - 2\Vert x - x^n \Vert
- \Vert x^n \Vert
$$
$$
\geq C \Vert h \Vert^p - 2\Vert x - x^n \Vert).
$$
Consider now $N$ large enough so that $\Vert x-x^N\Vert \le
\frac{C}{2}\delta^p$. For $h \in H_N$, $\Vert h\Vert=\delta$
we have
$$
\Vert x + h \Vert - \Vert x \Vert \geq
\frac{C}{2}\Vert h \Vert^p
$$
as we required.
\bigskip
In the second step we characterize the
spaces $\bigoplus\Sb \ell_{p} \endSb\ell_{q}^{n}$ which
admit a separating polynomial.

\bn
{\bf Step 2.} {\it
The space  $X=\bigoplus\Sb \ell_{p} \endSb\ell_{q}^{n}$
admits a separating polynomial if and only if $p,q$
are both even integers satisfying that $p=kq$,
for some $k \in \entier$.}
\bn
{\bf Proof of Step 2:}
First of all, it is clear that if $p$ is a multiple of $q$ and both are
even integers, the space $X$ has a separating polynomial. Indeed, the expression
$P(x) = \Vert x \Vert^p$ defines a $p$-homogeneous separating polynomial on
$X$.
\smallskip
On the other hand, if $X$ admits a separating polynomial, $p$ is an even
integer (since $\ell_p$ is a subspace of $X$).

We may assume that $p>q$. Otherwise, if $p<q$, since $X$ is saturated with
subspaces of cotype $p$, by [4] $X$ would actually have cotype $p$, but
$\ell_q $ is finitely represented in  $X$, and this is not possible.

Since $p>q$ the space $X$ has modulus of convexity of power
type $p$ (see [7]) and by Theorem 1-1,  there exists
a separating polynomial of degree at most $p$ on $X$. By Lemma
3-3, we may actually assume that there is a
$p$-homogeneous separating polynomial.

\smallskip
Let $P$ be a $p$-homogeneous polynomial on $X$
verifying that $P(x) \geq 1$ if $\| x \| \geq 1$.
Let us denote by $\pi_n$ the
projection of the space $\ell_q$ onto the finite dimensional subspace
$\ell_q^n$, and by $i_n$ the inclusion mapping from $\ell_q^n$ into the space
$X=\bigoplus_{\ell_p}\ell_{q}^{n}$.
We consider the sequence of polynomials $\{P_n\}_n$ on $\ell_q$
defined by $P_n=P\circ i_n\circ \pi_n$.
First of all, $\{P_n\}_n$ is a uniformly bounded sequence of $p$-homogeneous
polynomials on $X$; hence by Theorem 4 in [14] there is a subsequence
$\{P_{n_j}\}_j$  of  $\{P_n\}_n$ which converges pointwise to
a $p$-homogeneous polynomial $P^*$ on  $\ell_q$.
Let $x\in\ell_q$ be fixed,
with $\Vert x\Vert =2$; then for $n$ large enough, we have that
$P_n(x) \ge 1$.
This implies that the  $P^*$ is a $p$-homogeneous separating polynomial on
$\ell_{q}$, and  by [9] it follows that $p$ must be a multiple of
$q$.

\bn
{\bf Remark.} Since the space
$\bigoplus\Sb \ell_{2} \endSb\ell_{4}^{n}$ has
modulus of smoothness of power type 2,  by the results
in [12], it admits a twice  G\^ateaux differentiable equivalent norm
and, in particular, a continuous bump function with second-order directional
Taylor expansion at each
point. However, by Theorem 3-2, there is not even a bump function with
Taylor expansion of
second order at each point. In particular, this space can not be renormed with
an equivalent twice Fr\'echet differentiable norm.

\bn
{\bf Proof of Theorem 3.1:} First, in the case that  
$r=min\{p,q\}$ is not an even integer, since
$L_r \subset L_p(L_q)$ the result follows by Corollary 1-2.

Next, let  $r=q$ be an even integer. Since the space $X$ has modulus of
convexity of
power type $r$ (see [7]) by Theorem 1-1,  the existence of a $T^r$-smooth
bump function
implies the existence of a separating polynomial in $X$, and therefore
in the subspace $\bigoplus_{\ell_p}\ell_{q}^{n}$, which is not possible (
see Second
Step in the proof of Theorem 3-2).

Finally, in the case that $r=p$ is an even integer, if the space
$X$ has a $T^p$-smooth bump function, also has the subspace
$\bigoplus_{\ell_p}\ell_{q}^{n}$.
Then, it follows from Theorem 3-2 that $p,q$ must be both even integers,
and $p$ is a multiple of $q$.

\bigskip
In the following Proposition we prove that for each real valued function in
the space $X=\bigoplus_{\ell_p}\ell_{q}^{n}$ with certain properties of uniform
higher differentiability, we may contruct an associated function on the
space $\ell_q$
with shares the same properties of differentiability.

\bn
{\bf Proposition 3.4.} {\it
Let $1<p,q< \infty $ and  $X=\bigoplus_{\ell_p}\ell_{q}^{n}$. Let   $F$ be
a real valued function on  $X$  which is uniformly   $H^r$-smooth. For each
$n$ we consider the function  $F_n$ defined on  $\ell_q$ as  $F_n = F\circ
i_n\circ \pi_n$ .
Then there exists a subsequence  $\{F_{n_j}\}_j$  of  $\{F_n\}_n$ which
converges to
a uniformy $H^r$-smooth function  $F^*$ on  $\ell_q$. }

\bigskip\noindent
{\bf Proof : }
Since $F$ is uniformly  $H^r$-smooth, if $k$ is the greatest
integer strictly
less than $r$, $F$ is $C^k$-smooth and
there exists a constant $C>0$ such that for all
$x,y\in X$,
$$
\Vert F^{(k)} (x) -  F^{(k)}(y) \Vert \leq C\Vert x-y \Vert^{r-k}.
$$
We claim that the family  $\{ F^{(i)}_{n}\}$ for $ 0 \leq  i \leq k$, is
equicontinuous
on bounded sets of $\ell_q$. Indeed, denote $T_n = i_n \circ\pi_n$; then
for all $x,y,h \in \ell_q$, with $\Vert h \Vert \leq 1$,
$$
\big\vert \Bigl( F^{(k)}_n(x) -  F^{(k)}_n(y) \Bigr)(h) \big\vert =
\big\vert \Bigl( F^{(k)}_n\bigl(T_n(x)\bigr) -
F^{(k)}_n\bigl(T_n(y)\bigr)\Bigr)(h)\big\vert
\leq\qquad (*)
$$
$$
\leq C\Vert T_n(x) - T_n(y) \Vert^{r-k}.
\Vert T_n(h)\Vert  \leq
C \Vert x-y \Vert^{r-k}.
$$
\bigskip

Hence the sequence  $ \{ F^{(k)}_{n} \}_n$ is equicontinuous on
$\ell_q$, and also uniformly bounded on bounded sets of $\ell_q$.
This implies, by using the mean value theorem,
that the sequence  $ \{ F^{(i)}_{n} \}$ for $0\leq i \leq k$ is
equicontinuous on bounded sets of $\ell_q$.
Then, by using Lemma 4 in [14], for each $m \in \entier$ 
there exists a subsequence of $\{F_ n\}_n$ which converges 
pointwise to a function $F_m^*$ on the ball $B(0;m)$ of $\ell_q$, and
with their $i$-th derivatives pointwise converging to the $i$-th 
derivative of $F_m^*$ for all $i=0,\dots,k$. 
By a diagonalization procedure, passing to a subsequence,
we may assume that there exists a subsequence $\{F_{n_j}\}_j$ 
converging to a function $F^*$ pointwise on the whole space $\ell_q$. 
Then, by passing to the limit in (*), it follows that $F^*$ is also 
uniformly $H^r$-smooth, as we required.

\bigskip
From this proposition we obtain the following corollary.

\bn
{\bf
Corollary 3.5. }{\it Let  $1<q<p<\infty$, and $q$ not an even integer. Then
the space  $ X=\bigoplus_{\ell_p} \ell_{q}^{n}$  is not $H^r$-smooth for
$r>q$.}

\bigskip\noindent
{\bf Proof : } Assume that $X$ is $H^r$-smooth.
Since the space $X$ does not contain an isomorphic copy of $c_0$, it
is indeed uniformly  $H^r$-smooth (see [5, Th. V-3-1]). Let $F$ be
a function on $X$ which is uniformly  $H^r$-smooth, and such that
$F(x)=0$  whenever $\Vert x \Vert \geq 1$ and  $F(0)=1$.
Consider the sequence of functions  $\{F_ n\}_n$ as in Proposition 3-4.
Then, for a fixed point $x$ in $\ell_q$ with $\Vert x \Vert\geq 2$,
for $n$ large enough we have that
$F_n(x) =0$. This implies that the function $F^*$ obtained according to
Proposition 3-4 is, actually,  a uniformly  $H^r$-smooth bump function on 
$\ell_q$. But since $r>q$ this is not possible.
\bigskip
{\bf Remark.} From Theorem 3-2 and Corollary 3-5, we obtain
that, in Theorem 3-1, we may replace  $\bigoplus_{\ell_p} \ell_{q}^{n}$ by
$L_p(L_q )$, and the conclusion is the same except in the case that
$q<p$, and $q$ is not an even integer. In this case, we know from 
Corollary 3-5 that there is no uniformly $H^r$-smooth bump function 
for $r>q$. However, we do not know whether it is also true in this case that 
there is no $T^q$-smooth bump function.

\bigskip
\centerline{\bf References}

\medskip\noindent
[1] AVEZ A., {\it Calcul diff\'erentiel}. Masson, 1983.

\medskip\noindent
[2] BONIC R. and FRAMPTON J., Smooth functions on Banach manifolds,
{\it J. Math. Mech.} {\bf 15}, 877-898, (1966).

\medskip\noindent
[3] CARTAN H., {\it Cours de calcul diff\'erentiel}. Hermann (Collection
M\'ethodes), 1977.

\medskip\noindent
[4] DEVILLE R., Geometrical implications of the existence of very smooth
bump functions on Banach spaces,
{\it Israel. J. Math.} {\bf 67}, 1-22, (1989).

\medskip\noindent
[5] DEVILLE R., GODEFROY G. and ZIZLER V., {\it Smoothness and renormings
in Banach spaces}. Pitman Monographs in Mathematics 64. 
Longman Scientific and Technical, 1993.

\medskip\noindent
[6] FABIAN M., PREISS D., WHITFIELD J. H. M. and ZIZLER V.,
Separating polynomials on Banach spaces,
{\it Quart. J. Math. Oxford} {\bf 40}, 409-422, (1989).

\medskip\noindent
[7] FIGIEL T., On the moduli of convexity and smoothness,
{\it Studia Math.} {\bf 56}, 121-155, (1976).

\medskip\noindent
[8] GONZALO R. and JARAMILLO J.A., Compact polynomials between Banach spaces,
{\it Proc. Roy. Irish Acad.} {\bf 95A}, 213-226, (1995).

\medskip\noindent
[9] HABALA P. and HAJEK P., Stabilization of polynomials, {\it C. R. Acad
Sci. Paris} {\bf 320}, 821-825, (1995).

\medskip\noindent
[10] KURZWEIL J., On approximation in real Banach spaces, {\it Studia
Math.} {\bf 14}, 213-231, (1954).

\medskip\noindent
[11] LEONARD I.E. and SUNDARESAN K., Geometry of Lebesgue-Bochner function
spaces-Smoothness,
{\it Trans. Amer. Math. Soc.} {\bf 198}, 229-251, (1974).

\medskip\noindent
[12] McLAUGHLIN, R. POLIQUIN, J. VANDERWEFF,J. and
ZIZLER V., Second-order G\^ateaux differentiable bump functions and
approximations in Banach spaces, 
{\it  Can. J. Math.} {\bf 45(3)}, 612-625, (1993).

\medskip\noindent
[13] MALEEV, R. M. and TROYANSKI, S. L., Smooth functions in Orlicz spaces,
{\it Contemporary Math.} {\bf 85}, 355-369, (1989).

\medskip\noindent
[14] NEMIROVSKII A. S. and SEMENOV S. M.,
On polynomial approximation of functions on Hilbert spaces,
{\it Math. USSR Sbornik} {\bf 21}, 255-277, (1973).

\medskip\noindent
[15] PHELPS R. R., {\it Convex functions, monotone operators and
differentiability}. Lecture Notes in Mathematics 1364. 
Springer-Verlag, 1989.

\medskip\noindent
[16] STEGALL C., Optimization of functions on certain subsets of Banach spaces,
{\it Math. Ann.} {\bf 236}, 171-176, (1978).

\enddocument